\documentclass[preprint,12pt]{elsarticle}

\usepackage{amssymb}
\usepackage{amsthm}
\usepackage{graphics}
\usepackage{epsfig}
\usepackage{graphicx}

\usepackage{a4wide}
\usepackage{amsmath}
\usepackage{amsfonts}
\usepackage{enumerate}
\newcommand{\pf}{\noindent {\bf Proof: }}
\newtheorem*{theorem*}{Main Theorem}
\newtheorem{lemma}{\bf Lemma}[section]
 
\newtheorem{theorem}{\bf Theorem}[section]
\newtheorem{proposition}[lemma]{\bf Proposition}

\begin{document}

\begin{frontmatter}



\title{On Perfectness of Intersection Graph of Ideals of $\mathbb{Z}_n$}



\author{Angsuman Das\corref{cor1}}
\ead{angsumandas@sxccal.edu}

\address{Department of Mathematics,\\ St.Xavier's College, Kolkata, India.\\angsumandas@sxccal.edu}
\cortext[cor1]{Corresponding author}


\begin{abstract}
In this paper, we characterize the positive integers $n$ for which intersection graph of ideals of $\mathbb{Z}_n$ is perfect. 
\end{abstract}


\end{frontmatter}


\section{Introduction}
The idea of associating graphs to algebraic structures for characterizing the algebraic structures with graphs and vice versa dates back to Bosak \cite{bosak}. Till then, a lot of research, e.g., \cite{graph-ideal,anderson-livingston,badawi,power2,mks-ideal,angsu-comm-alg-1,angsu-lin-mult-alg,angsu-comm-alg-2,angsu-jaa,survey2} has been done in connecting graph structures to various algebraic objects like groups, rings, vector spaces etc. However, the most prominent among them are the zero-divisor graphs \cite{anderson-livingston} and intersection graph of ideals of rings  \cite{mks-ideal}.  Recently, authors in \cite{weakly-perfect} proved that intersection graph of ideals of $\mathbb{Z}_n$ is weakly perfect for all $n>0$. In this paper, we characterize the values of $n$ for which the intersection graph of ideals of $\mathbb{Z}_n$ is perfect. In particular, we prove the following theorem.

\begin{theorem*}
The intersection graph of ideals of $\mathbb{Z}_n$ is perfect if and only if  $n={p_1}^{\alpha_1}{p_2}^{\alpha_2}{p_3}^{\alpha_3}{p_4}^{\alpha_4}$ where $p_i$'s are distinct primes and $\alpha_i \in \mathbb{N}\cup\{0\}$, i.e., the number of distinct prime factors of $n$ is less than or equal to $4$.
\end{theorem*}

\section{Definition, Preliminaries and Known Results}
In this section, for convenience of the reader and also for later use, we recall some definitions, notations and results concerning elementary graph theory and intersection graph of ideals of a ring. For undefined terms and concepts the reader is referred to \cite{west-graph-book}.

By a graph $G=(V,E)$, we mean a non-empty set $V$ and a symmetric binary relation (possibly empty) $E$ on $V$. The set $V$ is called the set of vertices and $E$ is called the set of edges of $G$. Two element $u$ and $v$ in $V$ are said to be adjacent if $(u,v) \in E$. $H=(W,F)$ is called an {\it induced subgraph} of $G$ if $\phi \neq W \subseteq V$ and $F$ consists of all the edges between the vertices in $W$ in $G$.  A complete subgraph of a graph $G$ is called a {\it clique}. A {\it maximal clique} is a clique which is maximal with respect to inclusion. The {\it clique number} of $G$, written as $\omega(G)$, is the maximum size of a clique in $G$. The {\it chromatic number} of $G$, denoted as $\chi(G)$, is the minimum number of colours needed to label the vertices so that the adjacent vertices receive different colours. It is easy to observe that $\omega(G)\leq \chi(G)$. A graph $G$ is said to be {\it weakly perfect} if $\omega(G)= \chi(G)$ and it is said to be {\it perfect} if $\omega(H)= \chi(H)$ for all induced subgraphs $H$ of $G$. Chudnovsky {\it et.al.} \cite{spgt} in 2004 settled a long standing conjecture regarding perfect graphs and provided a characterization of perfect graphs.

{\theorem[Strong Perfect Graph Theorem] \label{spgt} \cite{spgt} A graph $G$ is perfect if and only if neither $G$ nor its complement contains an odd cycle of length at least $5$ as an induced subgraph.}

Let $R$ be a ring. The intersection graph of ideals of $R$ (introduced in \cite{mks-ideal}), denoted by $G(R)$, consists of all non-trivial ideals as vertices and two ideals $I$ and $J$ are adjacent if and only if $I \cap J\neq \{0\}$. Throughout this paper, we take the ring $R$ to be $\mathbb{Z}_n$, the ring of integers modulo $n$. We know that $\mathbb{Z}_n$ is a principal ideal ring and each of its ideals is
generated by $\overline{m}\in \mathbb{Z}_n$ where $m$ is a factor of $n$. For convenience, we denote this ideal by $(m)$. Also without loss of generality, whenever we take an ideal $(m)$ of $\mathbb{Z}_n$, we assume that $m$ is a factor of $n$. It was proved in \cite{weakly-perfect} proved that intersection graph of ideals of $\mathbb{Z}_n$ is weakly perfect, i.e., $\omega(G( \mathbb{Z}_n))=\chi(G( \mathbb{Z}_n))$ for all $n>0$.

\section{Perfectness of Intersection Graph of Ideals of $\mathbb{Z}_n$}
In this section, we prove some preparatory results and subsequently use them to prove the main theorem of the paper.
{\proposition \label{lcm-lemma} Let $G(\mathbb{Z}_n)$ be the intersection graph of ideals of $\mathbb{Z}_n$ and $(a)$ and $(b)$ be two ideals in $\mathbb{Z}_n$ such that $a \mid n$ and $b\mid n$. Then $(a)$ and $(b)$ are adjacent in $G(\mathbb{Z}_n)$  if and only if $lcm(a,b)$ is a factor of $n$ and $1<lcm(a,b)<n$.}\\
\\
\pf Since $\mathbb{Z}_n$ is isomorphic to $\mathbb{Z}/n\mathbb{Z}$ as ring via the correspondence $\overline{a}\leftrightarrow a+n\mathbb{Z}$, the ideal $(a)$ in $\mathbb{Z}_n$ corresponds to the ideal $\langle a \rangle + n\mathbb{Z}$ in $\mathbb{Z}/n\mathbb{Z}$ where $\langle a \rangle$ denote the set of integer multiples of $a$. Now, let $(a)\sim (b)$ in $G(\mathbb{Z}_n)$, i.e., $(a)\cap(b)\neq \{\overline{0}\}$. Since $a \mid n$ and $b\mid n$, we have $lcm(a,b)\mid n$. On the other hand, using the correspondence described above, we have  $\langle a \rangle + n\mathbb{Z} \cap \langle b \rangle + n\mathbb{Z}\neq \{n\mathbb{Z}\}$. But, we know that $\langle a \rangle + n\mathbb{Z} \cap \langle b \rangle + n\mathbb{Z}=\langle lcm(a,b) \rangle + n\mathbb{Z}$. Hence, we have $\langle lcm(a,b) \rangle + n\mathbb{Z} \neq \{n\mathbb{Z}\}$. This, together with the fact that $lcm(a,b)\mid n$, implies that $1<lcm(a,b)<n$.

Conversely, let $lcm(a,b)$ is a factor of $n$ and $1<lcm(a,b)<n$. Clearly, $\overline{0} \neq \overline{lcm(a,b)} \in (a) \cap (b)$ in $\mathbb{Z}_n$ and hence $(a)\sim (b)$ in $G(\mathbb{Z}_n)$.  \qed

{\theorem \label{not-perfect} Let $n={p_1}^{\alpha_1}{p_2}^{\alpha_2}\cdots {p_k}^{\alpha_k}$. If $k\geq 5$, then $G(\mathbb{Z}_n)$ is not perfect.}\\
\\
\pf Let $n={p_1}^{\alpha_1}{p_2}^{\alpha_2}\cdots {p_5}^{\alpha_5}.s$ where $s=1$ if $k=5$ and $s={p_6}^{\alpha_6}\cdots {p_k}^{\alpha_k}$ if $k>5$. Consider the cycle $C$ given by $({p_1}^{\alpha_1}{p_2}^{\alpha_2} {p_3}^{\alpha_3}s)\sim({p_2}^{\alpha_2}{p_3}^{\alpha_3} {p_4}^{\alpha_4}s)\sim({p_3}^{\alpha_3}{p_4}^{\alpha_4} {p_5}^{\alpha_5}s)\sim({p_4}^{\alpha_4}{p_5}^{\alpha_5} {p_1}^{\alpha_1}s)\sim({p_5}^{\alpha_5}{p_1}^{\alpha_1}{p_2}^{\alpha_2} s)\sim ({p_1}^{\alpha_1}{p_2}^{\alpha_2} {p_3}^{\alpha_3}s)$. Simple calculation using Proposition \ref{lcm-lemma} shows that $C$ is an induced $5$-cycle in $G(\mathbb{Z}_n)$ and hence by Theorem \ref{spgt}, $G(\mathbb{Z}_n)$ is not perfect.\qed

{\theorem \label{no-odd-hole} Let $n={p_1}^{\alpha_1}{p_2}^{\alpha_2}{p_3}^{\alpha_3}{p_4}^{\alpha_4}$. Then $G(\mathbb{Z}_n)$ does not contain any induced cycle of length greater than $4$.}\\
\\
\pf Let, if possible $G(\mathbb{Z}_n)$  contains an induced cycle $C$ of length greater than $4$, say $(a_1)\sim(a_2)\sim(a_3)\sim(a_4)\sim(a_5)\sim \cdots \sim (a_1)$. By Proposition \ref{lcm-lemma}, we have $$lcm(a_1,a_3)=lcm(a_1,a_4)=lcm(a_2,a_4)=lcm(a_2,a_5)=lcm(a_3,a_5)=n.$$

{\it Claim: $gcd(a_1,a_3)>1$.} If possible, let $gcd(a_1,a_3)=1$. Since $gcd(a_1,a_3)\cdot lcm(a_1,a_3)=a_1a_3$, we have $lcm(a_1,a_3)=a_1a_3=n$. Note that as $lcm(a_3,a_5)=n$, we have $n=\frac{a_1a_3}{gcd(a_1,a_3)}=\frac{a_3a_5}{gcd(a_3,a_5)}$, i.e., $a_1\cdot gcd(a_3,a_5)=a_5\cdot gcd(a_1,a_3)$, i.e., $a_5=a_1\cdot gcd(a_3,a_5)$, i.e., $a_5$ is a multiple of $a_1$. Now as $a_1$ and $a_3$ are coprime and their lcm is $n$, without loss of generality, two cases may arise: either $a_1={p_1}^{\alpha_1}{p_2}^{\alpha_2}; a_3={p_3}^{\alpha_3}{p_4}^{\alpha_4}$ or $a_1={p_1}^{\alpha_1};a_3={p_2}^{\alpha_2}{p_3}^{\alpha_3}{p_4}^{\alpha_4}$.

If $a_1={p_1}^{\alpha_1}{p_2}^{\alpha_2}; a_3={p_3}^{\alpha_3}{p_4}^{\alpha_4}$, we have $a_5={p_1}^{\alpha_1}{p_2}^{\alpha_2}\cdot s$ for some natural number $s$ such that $a_5 \mid n$. Also as $lcm(a_1,a_4)=n$, we have $a_4={p_3}^{\alpha_3}{p_4}^{\alpha_4}\cdot t$, for some natural number $t$ such that $a_4 \mid n$. Thus $lcm(a_4,a_5)=n$ contradicting Proposition \ref{lcm-lemma} and the fact that $a_4\sim a_5$ in $C$. 

If $a_1={p_1}^{\alpha_1};a_3={p_2}^{\alpha_2}{p_3}^{\alpha_3}{p_4}^{\alpha_4}$, similarly we have $a_5={p_1}^{\alpha_1}\cdot s$ and $a_4={p_2}^{\alpha_2}{p_3}^{\alpha_3}{p_4}^{\alpha_4} \cdot t$ and hence $lcm(a_4,a_5)=n$ thereby leading to a contradiction. Thus by combining above two cases, we have $gcd(a_1,a_3)>1$. 

Thus we have $lcm(a_1,a_3)=n$ and $gcd(a_1,a_3)>1$ with  $a_1 \mid n$ and $a_3\mid n$. Without loss of generality, let $p_1$ be a common factor of $a_1$ and $a_3$ and let $a_1={p_1}^x\cdot s$ and $a_3={p_1}^y \cdot t$ where $p_1$ is coprime with $s$ and $t$. Now, if $max\{x,y\}<\alpha_1$, then $lcm(a_1,a_3)<n$, a contradiction. Thus either $x=\alpha_1$ or $y=\alpha_1$, i.e., for any common prime divisor $p_i$ of $a_1$ and $a_3$, either ${p_i}^{\alpha_i} \mid a_1$ or ${p_i}^{\alpha_i} \mid a_3$ or both. Also as $lcm(a_1,a_3)=n$, all the ${p_i}^{\alpha_i}$ are factors of either $a_1$ or $a_3$ or both. Thus, without loss of generality, the forms of $a_1$ and $a_3$ are as follows: either $$\mathsf{Case ~ 1:}~ a_1={p_1}^{\alpha_1}{p_2}^{\alpha_2}{p_3}^{\beta_3}{p_4}^{\beta_4};a_3={p_1}^{\beta_1}{p_2}^{\beta_2}{p_3}^{\alpha_3}{p_4}^{\alpha_4}$$  or $$\mathsf{Case ~ 2:}~a_1={p_1}^{\alpha_1}{p_2}^{\beta_2}{p_3}^{\beta_3}{p_4}^{\beta_4};a_3={p_1}^{\beta_1}{p_2}^{\alpha_2}{p_3}^{\alpha_3}{p_4}^{\alpha_4}$$ or $$\mathsf{Case ~ 3:}~a_1={p_1}^{\alpha_1}{p_2}^{\alpha_2}{p_3}^{\beta_3}{p_4}^{\beta_4};a_3={p_1}^{\beta_1}{p_2}^{\alpha_2}{p_3}^{\alpha_3}{p_4}^{\alpha_4}$$ or $$\mathsf{Case ~ 4:}~a_1={p_1}^{\alpha_1}{p_2}^{\alpha_2}{p_3}^{\alpha_3}{p_4}^{\beta_4};a_3={p_1}^{\beta_1}{p_2}^{\alpha_2}{p_3}^{\alpha_3}{p_4}^{\alpha_4}$$ where $\beta_i < \alpha_i$. Note that in first two cases, $a_1$ and $a_3$ do not share any ${p_i}^{\alpha_i}$ as common factor. In the third case, they share only one ${p_i}^{\alpha_i}$ as common factor and in the fourth case, they share two ${p_i}^{\alpha_i}$'s as common factor.

{\it Case 1: ($a_1={p_1}^{\alpha_1}{p_2}^{\alpha_2}{p_3}^{\beta_3}{p_4}^{\beta_4};a_3={p_1}^{\beta_1}{p_2}^{\beta_2}{p_3}^{\alpha_3}{p_4}^{\alpha_4}$)} Since $lcm(a_1,a_4)=n$, we have $a_4={p_1}^{\gamma_1}{p_2}^{\gamma_2}{p_3}^{\alpha_3}{p_4}^{\alpha_4}$ where $\gamma_1\leq \alpha_1,\gamma_2 \leq \alpha_2$ and $(\gamma_1,\gamma_2)\neq (\alpha_1,\alpha_2)$. Again, since $lcm(a_3,a_5)=n$, we have $a_5={p_1}^{\alpha_1}{p_2}^{\alpha_2}{p_3}^{\delta_3}{p_4}^{\delta_4}$ where $\delta_3\leq \alpha_3,\delta_4 \leq \alpha_4$ and $(\delta_3,\delta_4)\neq (\alpha_3,\alpha_4)$. Hence, we have $lcm(a_4,a_5)=n$, a contradiction to the fact that $a_4 \sim a_5$. Thus Case 1 is an impossibility.

{\it Case 2: ($a_1={p_1}^{\alpha_1}{p_2}^{\beta_2}{p_3}^{\beta_3}{p_4}^{\beta_4};a_3={p_1}^{\beta_1}{p_2}^{\alpha_2}{p_3}^{\alpha_3}{p_4}^{\alpha_4}$)} Since $lcm(a_1,a_4)=n$, we have $a_4={p_1}^{\gamma_1}{p_2}^{\alpha_2}{p_3}^{\alpha_3}{p_4}^{\alpha_4}$ where $\gamma_1< \alpha_1$. Again, since $lcm(a_3,a_5)=n$, we have $a_5={p_1}^{\alpha_1}{p_2}^{\delta_2}{p_3}^{\delta_3}{p_4}^{\delta_4}$ where $\delta_i\leq \alpha_i$ and $(\delta_2,\delta_3,\delta_4)\neq (\alpha_2,\alpha_3,\alpha_4)$. Hence, we have $lcm(a_4,a_5)=n$, a contradiction to the fact that $a_4 \sim a_5$. Thus Case 2 is an impossibility.

{\it Case 3: ($a_1={p_1}^{\alpha_1}{p_2}^{\alpha_2}{p_3}^{\beta_3}{p_4}^{\beta_4};a_3={p_1}^{\beta_1}{p_2}^{\alpha_2}{p_3}^{\alpha_3}{p_4}^{\alpha_4}$)} Since $lcm(a_1,a_4)=n$, we have ${p_3}^{\alpha_3}{p_4}^{\alpha_4}\mid a_4$. Again, since $lcm(a_3,a_5)=n$, we have ${p_1}^{\alpha_1}\mid a_5$. Now, as $lcm(a_2,a_5)=n$, we have either ${p_2}^{\alpha_2}\mid a_2$ or ${p_2}^{\alpha_2}\mid a_5$. But if ${p_2}^{\alpha_2}\mid a_5$, then we have $lcm(a_4,a_5)=n$, a contradiction. Thus, we have ${p_2}^{\alpha_2}\mid a_2$. Again, as $lcm(a_2,a_4)=n$, we have either ${p_1}^{\alpha_1}\mid a_2$ or ${p_1}^{\alpha_1}\mid a_4$. If ${p_1}^{\alpha_1}\mid a_2$, then $lcm(a_2,a_3)=n$, a contradiction. On the other hand, if ${p_1}^{\alpha_1}\mid a_4$, then $lcm(a_3,a_4)=n$, a contradiction. Thus Case 3 is an impossibility.

{\it Case 4: ($a_1={p_1}^{\alpha_1}{p_2}^{\alpha_2}{p_3}^{\alpha_3}{p_4}^{\beta_4};a_3={p_1}^{\beta_1}{p_2}^{\alpha_2}{p_3}^{\alpha_3}{p_4}^{\alpha_4}$)} Since $lcm(a_1,a_4)=n$, we have ${p_4}^{\alpha_4}\mid a_4$.  Now, as $lcm(a_2,a_4)=n$, we have either ${p_1}^{\alpha_1}\mid a_2$ or ${p_1}^{\alpha_1}\mid a_4$. If ${p_1}^{\alpha_1}\mid a_2$, then $lcm(a_2,a_3)=n$, a contradiction. On the other hand, if ${p_1}^{\alpha_1}\mid a_4$, then $lcm(a_3,a_4)=n$, a contradiction. Thus Case 4 is an impossibility.

Thus, combining all the cases we conclude that $G(\mathbb{Z}_n)$ does not contain any induced cycle of length greater than $4$.\qed

{\theorem \label{no-odd-antihole} Let $n={p_1}^{\alpha_1}{p_2}^{\alpha_2}{p_3}^{\alpha_3}{p_4}^{\alpha_4}$. Then $\overline{G(\mathbb{Z}_n)}$, the complement of $G(\mathbb{Z}_n)$, does not contain any induced cycle of length greater than $4$.}\\
\\ \pf Let, if possible $\overline{G(\mathbb{Z}_n)}$  contains an induced cycle $C$ of length greater than $4$, say $(a_1)\sim(a_2)\sim(a_3)\sim(a_4)\sim \cdots \sim(a_t)\sim (a_1)$ with $t\geq 5$. Then, by Proposition \ref{lcm-lemma}, $lcm(a_1,a_2)=lcm(a_2,a_3)=lcm(a_3,a_4)=\cdots =lcm(a_t,a_1)=n$.

[{\it Claim: $gcd(a_2,a_3)>1$}] If possible, let $gcd(a_2,a_3)=1$. Since $lcm(a_2,a_3)=n$, we have $n=a_2a_3$. Thus without loss of generality, either $$a_2={p_1}^{\alpha_1}{p_2}^{\alpha_2};a_3={p_3}^{\alpha_3}{p_4}^{\alpha_4} \mbox{ or } a_2={p_1}^{\alpha_1};a_3={p_2}^{\alpha_2}{p_3}^{\alpha_3}{p_4}^{\alpha_4}$$

If $a_2={p_1}^{\alpha_1}{p_2}^{\alpha_2};a_3={p_3}^{\alpha_3}{p_4}^{\alpha_4}$, as $lcm(a_3,a_4)=lcm(a_1,a_2)=n$, we have $a_1={p_3}^{\alpha_3}{p_4}^{\alpha_4}\cdot s$ and $a_4={p_1}^{\alpha_1}{p_2}^{\alpha_2}\cdot t$ for some positive integer $s,t$. But this implies that $lcm(a_1,a_4)=n$, i.e., $a_1\sim a_4$ in $\overline{G(\mathbb{Z}_n)}$, a contradiction.

On the other hand, if $a_2={p_1}^{\alpha_1};a_3={p_2}^{\alpha_2}{p_3}^{\alpha_3}{p_4}^{\alpha_4}$, as $lcm(a_3,a_4)=lcm(a_1,a_2)=n$, we have $a_1={p_2}^{\alpha_2}{p_3}^{\alpha_3}{p_4}^{\alpha_4}\cdot s$ and $a_4={p_1}^{\alpha_1}\cdot t$ for some positive integer $s,t$. But this implies that $lcm(a_1,a_4)=n$, i.e., $a_1\sim a_4$ in $\overline{G(\mathbb{Z}_n)}$, a contradiction. Hence the claim is true.

Now, we have $lcm(a_2,a_3)=n$ and $gcd(a_2,a_3)>1$ with  $a_2 \mid n$ and $a_3\mid n$. Without loss of generality, let $p_1$ be a common factor of $a_2$ and $a_3$ and let $a_2={p_1}^x\cdot s$ and $a_3={p_1}^y \cdot t$ where $p_1$ is coprime with $s$ and $t$. Now, if $max\{x,y\}<\alpha_1$, then $lcm(a_2,a_3)<n$, a contradiction. Thus either $x=\alpha_1$ or $y=\alpha_1$, i.e., for any common prime divisor $p_i$ of $a_2$ or $a_3$, either ${p_i}^{\alpha_i} \mid a_2$ or ${p_i}^{\alpha_i} \mid a_3$ or both. Also as $lcm(a_2,a_3)=n$, all the ${p_i}^{\alpha_i}$ are factors of either $a_2$ or $a_3$. Thus, without loss of generality, the forms of $a_2$ and $a_3$ are as follows: either $$\mathsf{Case ~ 1:}~ a_2={p_1}^{\alpha_1}{p_2}^{\alpha_2}{p_3}^{\beta_3}{p_4}^{\beta_4};a_3={p_1}^{\beta_1}{p_2}^{\beta_2}{p_3}^{\alpha_3}{p_4}^{\alpha_4}$$  or $$\mathsf{Case ~ 2:}~a_2={p_1}^{\alpha_1}{p_2}^{\beta_2}{p_3}^{\beta_3}{p_4}^{\beta_4};a_3={p_1}^{\beta_1}{p_2}^{\alpha_2}{p_3}^{\alpha_3}{p_4}^{\alpha_4}$$ or $$\mathsf{Case ~ 3:}~a_2={p_1}^{\alpha_1}{p_2}^{\alpha_2}{p_3}^{\beta_3}{p_4}^{\beta_4};a_3={p_1}^{\beta_1}{p_2}^{\alpha_2}{p_3}^{\alpha_3}{p_4}^{\alpha_4}$$ or $$\mathsf{Case ~ 4:}~a_2={p_1}^{\alpha_1}{p_2}^{\alpha_2}{p_3}^{\alpha_3}{p_4}^{\beta_4};a_3={p_1}^{\beta_1}{p_2}^{\alpha_2}{p_3}^{\alpha_3}{p_4}^{\alpha_4}$$ where $\beta_i < \alpha_i$. Note that in first two cases, $a_2$ and $a_3$ do not share any ${p_i}^{\alpha_i}$ as common factor. In the third case, they share only one ${p_i}^{\alpha_i}$ as common factor and in the fourth case, they share two ${p_i}^{\alpha_i}$'s as common factor.

{\it Case 1: ($a_2={p_1}^{\alpha_1}{p_2}^{\alpha_2}{p_3}^{\beta_3}{p_4}^{\beta_4};a_3={p_1}^{\beta_1}{p_2}^{\beta_2}{p_3}^{\alpha_3}{p_4}^{\alpha_4}$)} Since $lcm(a_1,a_2)=lcm(a_3,a_4)=n$, we have ${p_3}^{\alpha_3}{p_4}^{\alpha_4} \mid a_1$ and 
 ${p_1}^{\alpha_1}{p_2}^{\alpha_2}\mid a_4$. But this implies $lcm(a_1,a_4)=n$, i.e., $a_1\sim a_4$ in $\overline{G(\mathbb{Z}_n)}$, a contradiction and hence Case 1 is an impossibility.
 
{\it Case 2: ($a_2={p_1}^{\alpha_1}{p_2}^{\beta_2}{p_3}^{\beta_3}{p_4}^{\beta_4};a_3={p_1}^{\beta_1}{p_2}^{\alpha_2}{p_3}^{\alpha_3}{p_4}^{\alpha_4}$)} Since $lcm(a_1,a_2)=lcm(a_3,a_4)=n$, we have ${p_2}^{\alpha_2}{p_3}^{\alpha_3}{p_4}^{\alpha_4}\mid a_1$ and ${p_1}^{\alpha_1}\mid a_4$. But this implies $lcm(a_1,a_4)=n$, a contradiction and hence Case 2 is an impossibility.
 
{\it Case 3: ($a_2={p_1}^{\alpha_1}{p_2}^{\alpha_2}{p_3}^{\beta_3}{p_4}^{\beta_4};a_3={p_1}^{\beta_1}{p_2}^{\alpha_2}{p_3}^{\alpha_3}{p_4}^{\alpha_4}$)}   Since $lcm(a_1,a_2)=n$, we have ${p_3}^{\alpha_3}{p_4}^{\alpha_4} \mid a_1$. Also, since $lcm(a_t,a_1)=n$, either ${p_1}^{\alpha_1}\mid a_1$ or ${p_1}^{\alpha_1}\mid a_t$. If ${p_1}^{\alpha_1}\mid a_1$, then we have ${p_1}^{\alpha_1}{p_3}^{\alpha_3}{p_4}^{\alpha_4} \mid a_1$ which implies $lcm(a_1,a_3)=n$, i.e., $a_1\sim a_3$ in $\overline{G(\mathbb{Z}_n)}$, a contradiction. On the other hand, if ${p_1}^{\alpha_1}\mid a_t$, we have $lcm(a_t,a_3)=n$, i.e., $a_t \sim a_3$ in $\overline{G(\mathbb{Z}_n)}$, a contradiction. Thus combining both the possibilities, Case 3 is an impossibility.

{\it Case 4: ($a_2={p_1}^{\alpha_1}{p_2}^{\alpha_2}{p_3}^{\alpha_3}{p_4}^{\beta_4};a_3={p_1}^{\beta_1}{p_2}^{\alpha_2}{p_3}^{\alpha_3}{p_4}^{\alpha_4}$)}   Since $lcm(a_1,a_2)=n$, we have ${p_4}^{\alpha_4} \mid a_1$. Also, since $lcm(a_t,a_1)=n$, either ${p_1}^{\alpha_1}\mid a_1$ or ${p_1}^{\alpha_1}\mid a_t$. If ${p_1}^{\alpha_1}\mid a_1$, then we have ${p_1}^{\alpha_1}{p_4}^{\alpha_4} \mid a_1$ which implies $lcm(a_1,a_3)=n$, i.e., $a_1\sim a_3$ in $\overline{G(\mathbb{Z}_n)}$, a contradiction. On the other hand, if ${p_1}^{\alpha_1}\mid a_t$, we have $lcm(a_t,a_3)=n$, i.e., $a_t \sim a_3$ in $\overline{G(\mathbb{Z}_n)}$, a contradiction. Thus combining both the possibilities, Case 4 is an impossibility.

Thus, combining all the cases we conclude that $\overline{G(\mathbb{Z}_n)}$ does not contain any induced cycle of length greater than $4$.\qed

Finally, with Theorems \ref{spgt}, \ref{not-perfect}, \ref{no-odd-hole} and \ref{no-odd-antihole} in hand, we are now in a position to prove the main result of this paper.
\begin{theorem*}
The intersection graph of ideals of $\mathbb{Z}_n$ is perfect if and only if  $n={p_1}^{\alpha_1}{p_2}^{\alpha_2}{p_3}^{\alpha_3}{p_4}^{\alpha_4}$ where $p_i$'s are distinct primes and $\alpha_i \in \mathbb{N}\cup\{0\}$, i.e., the number of distinct prime factors of $n$ is less than or equal to $4$.
\end{theorem*}
\pf Clearly, Theorem \ref{not-perfect} shows that the condition is necessary. For the sufficiency part, first with the help of Theorems \ref{no-odd-hole} and \ref{no-odd-antihole}, along with Theorem \ref{spgt}, we conclude that the intersection graph of ideals of $\mathbb{Z}_n$ is perfect if $n$ has exactly four distinct prime factors. The proofs for the cases when $n$ has exactly three, two or one distinct prime factors follows similarly by suitably taking some of the $\alpha_i$'s to be zero. \qed

\section*{Acknowledgement}
The author is thankful to Sabyasachi Dutta and Jyotirmoy Pramanik for some fruitful discussions on the paper. The research is partially funded by NBHM Research Project Grant, (Sanction No. 2/48(10)/2013/ NBHM(R.P.)/R\&D II/695), Govt. of India.

\end{document}